\documentclass[12pt]{amsart}
\usepackage{amsfonts}
\usepackage{amssymb}
\theoremstyle{definition}

\theoremstyle{remark}

\newcommand{\ds}{\displaystyle}

\begin{document}

\centerline{\large\bf ON THE BOCHNER CURVATURE TENSOR}
\centerline{\large\bf IN AN ALMOST HERMITIAN MANIFOLD
\footnote{\it SERDICA Bulgaricae mathematicae publicationes. Vol. 9, 1983, p. 168-171.}}

\vspace{0.2in}
\centerline{\large OGNIAN T. KASSABOV}

\vspace{0.5in}
{\sl We prove a classification theorem for $RK$-manifolds with linear dependence
between invariants of an antiholomorphic plane in the tangent space. As a 
consequence we find a characteristic condition for an $RK$-manifold to be of pointwise
constant antiholomorphic sectional curvature.}

\vspace{0.2in}
{\bf 1. Introduction.} Let $M$ be a $2n$-dimensional almost Hermitian manifolds, $n \ge 3$,
with metric tensor $g$ and almost complex structure $J$ and let $\nabla$ be the covariant
differentiation on $M$. The curvature tensor $R$ is defined by
$$
    R(X,Y,Z,U) = g(\nabla_X\nabla_Y Z - \nabla_Y\nabla_X Z - \nabla_{[X,Y]} Z,U)
$$
for $X,\, Y,\, Z,\, U \in \mathfrak X(M)$. The manifold is said to be an RK-manifold \cite{V} if
$$
	R(X,Y,Z,U) = R(JX,JY,JZ,JU)
$$
for all $X,Y,Z,U \in \mathfrak X(M)$. In this paper we trate for simplicity only the case 
of an RK-manifold although one cane made analogous considerations for an arbitrary almost
Hermitian manifold.

Let $ E_i ,\, i=1,\hdots,2n$ be a local orthonormal frame field. The Ricci tensor $S$
and the scalar curvature $\tau(R)$ are defined by
$$
    S(X,Y) = \sum_{i=1}^{2n} R(X,E_i,E_i,Y) \ , \quad \tau(R) = \sum_{i=1}^{2n} S(E_i,E_i).
$$
Analogously we set
$$
    S'(X,Y) = \sum_{i=1}^{2n} R(X,E_i,JE_i,JY) \ , \quad \tau'(R) = \sum_{i=1}^{2n} S'(E_i,E_i).
$$
We note that $S$ and $S'$ are symmetric and $S(X,Y)=S(JX,JY)$, $S'(X,Y)=S'(JX,JY)$.

The Bochner curvature tensor $B$ \cite{TV} for $M$ is defined by
$$
    B=R-\frac 1{8(n+2)}(\varphi + \psi)(S+3S')-\frac 1{8(n-2)}(3\varphi - \psi)(S-S')
$$
$$
    +\frac{\tau(R)+3\tau'(R)}{16(n+1)(n+2)}(\pi_1+\pi_2)
    +\frac{\tau(R)-\tau'(R)}{16(n-1)(n-2)}(3\pi_1-\pi_2),
$$
where $\varphi,\ \psi,\, \pi_1$ and $\pi_2$ are defined by 
$$
    \varphi(Q)(X,Y,Z,U) = g(X,U)Q(Y,Z) - g(X,Z)Q(Y,U)
$$
$$
    + g(Y,Z)Q(X,U) - g(Y,U)Q(X,Z)\, ,
$$
$$
    \psi(Q)(X,Y,Z,U) = g(X,JU)Q(Y,JZ) - g(X,JZ)Q(Y,JU)
$$
$$
     + g(Y,JZ)Q(X,JU) - g(Y,JU)Q(X,JZ)\, ,
$$
$$
    -2g(X,JY)Q(Z,JU) -2g(Z,JU)Q(X,JY)\, ,
$$
$$
    \pi_1(X,Y,Z,U) = g(X,U)g(Y,Z) - g(X,Z)g(Y,U)\, ,
$$
$$
    \pi_2(X,Y,Z,U) = g(X,JU)g(Y,JZ) - g(X,JZ)g(Y,JU) -2g(X,JY)g(Z,JU)\, .
$$

By a plane we mean a 2-dimensional linear subspace of the tangent space $T_p(M)$
of $M$ in $p$. A plane $\alpha$ is said to be holomorphic (resp. antiholomorphic)
if $J\alpha = \alpha$ (resp. $J\alpha$ is perpendicular to $\alpha$). 

A tensor field $T$ of type (0,4) is said to be an LC-tensor if it has the properties:

1) $ T(X,Y,Z,U) = - T(Y,X,Z,U)$,

2) $ T(X,Y,Z,U) = - T(X,Y,U,Z)$,

3) $ T(X,Y,Z,U) + T(Y,Z,X,U) + T(Z,X,Y,U) = 0$.

We need the following lemma.

L e m m a \cite{G}. {\it Let $M$ be a $2n$-dimensional almost Hermitian manifold, $n \ge 2$.
Let $T$ be an $LC$-tensor, satisfying the conditions:

1) $ T(X,Y,Z,U) = T(JX,JY,JZ,JU)$,

2) $ T(x,y,y,x) = 0$, where $\{ x,y \}$ is a basis of any holomorphic or antiholomorphic plane.

Then $T=0$.}

In section 2 we shall prove the following theorem.

T h e o r e m. {\it Let $M$ be a \ $2n$-dimensional RK-manifold, $n\ge 3$, which satisfies
$$
	\lambda R(x,y,y,x) + \mu (S(x,x)+S(y,y)) + \nu (S'(x,x)+S'(y,y)) = c(p) \leqno (1.1)
$$
for each point $p \in M$ and for all unit vectors $x,y \in T_p(M)$ with $g(x,y)=g(x,Jy)=0$,
where $\lambda,\, \mu,\, \nu$ are constants, $(\lambda,\mu,\nu) \ne (0,0,0)$ \ and \ $c(p)$
does not depend on \ $x,y$. Then

1) if $\lambda = 0$, then
$$
	\mu S+\nu S'= \frac{\mu \tau(R) + \nu \tau'(R)}{2n} g\ ;
$$

2) if  \ $\lambda \ne 0$, then \ $M$ \ has vanishing Bochner curvature tensor and the tensor \\
$((n+1)\lambda + 2(n^2-4)\mu)S + (2(n^2-4)\nu - 3\lambda)S'$ \ is proportional to the metric tensor}:
$$
	((n+1)\lambda + 2(n^2-4)\mu)S + (2(n^2-4)\nu - 3\lambda)S'
$$ 
$$
	=\frac{1}{2n}\{ ((n+1)\lambda + 2(n^2-4)\mu)\tau(R) + (2(n^2-4)\nu - 3\lambda)\tau'(R) \} g \ .
$$

An almost Hermitian manifold $M$ is said to be of pointwise constant antiholomorphic sectional
curvature if for each point $p \in M$ the curvature of an arbitrary antiholomorphic plane
$\alpha$ in $T_p(M)$ does not depend on $\alpha$.

C o r o l l a r y. {\it Let $M$ be a $2n$-dimensional RK-manifold, $n\ge 3$. Then $M$ has 
pointwise constant antiholomorphic sectional curvature if and only if $M$ has vanishing
Bochner curvature tensor and}
$$
	(n+1)S-3S'= \frac{1}{2n}((n+1)\tau(R) - 3\tau'(R))g \ .
$$

This is an analogue of a well known theorem of \, S c h o u t e n \, and \, S t r u i k \cite{S},
see also \cite{T}.

\vspace{.2in}
{\bf 2. Proof of the theorem.} Let \, $ e_i ,\, Je_i,\,  i=1,\hdots,2n$ \, be an orthonormal
basis of \, $T_p(M), \, p \in M$. In (1.1) we put \, $X=e_1, \, Y=e_i$ or $Y=Je_i$, $i=2,...,n$. 
Adding on $i$ we find
$$
	\lambda R(e_1,Je_1,Je_1,e_1) - (\lambda +2(n-2)\mu)S(e_1,e_1) - 2(n-2)\nu S'(e_1,e_1)
$$    
$$
	=\mu\tau(R)(p) + \nu\tau'(R)(p) - 2(n-1)c(p)
$$ 
and since we can take for $e_1$ an arbitrary unit vector in $T_p(M)$ we have
$$
	\begin{array}{c}
		\lambda H(x) - (\lambda +2(n-2)\mu)S(x,x) - 2(n-2)\nu S'(x,x) \\
		=\mu\tau(R)(p) + \nu\tau'(R)(p) - 2(n-1)c(p)
	\end{array}    \leqno(2.1)
$$ 
for each unit vector $x \in T_p(M)$, where $H(x)$ is the curvature of the holomorphic 
plane spanned by $x,\, Jx$, i.e. $H(x)=R(x,Jx,Jx,x)$.

If $\lambda = 0$ (2.1) takes the form
$$
	\mu S(x,x) + \nu S'(x,x) 
		=\frac{2(n-1)c(p)-\mu\tau(R)(p) - \nu\tau'(R)(p)}{2(n-2)} \ .
$$

We put $ x=e_i,\, x=Je_i$ and adding on $i$ we obtain
$$
	c=\frac{\mu\tau(R) + \nu\tau'(R)}{n}
$$
and case 1) is proved.

If $\lambda \ne 0$, (1.1) and (2.1) take the form
$$
	 R(x,y,y,x) + \mu_1 (S(x,x)+S(y,y)) + \nu_1 (S'(x,x)+S'(y,y)) = c_1(p) \, ,\leqno (2.2)
$$
$$
	\begin{array}{c}
		H(x) - (1 +2(n-2)\mu_1)S(x,x) - 2(n-2)\nu_1 S'(x,x) \\
		=\mu_1\tau(R)(p) + \nu_1\tau'(R)(p) - 2(n-1)c_1(p) \, ,
	\end{array}    \leqno(2.3)
$$ 
where $\mu_1=\mu/\lambda$, $\nu_1=\nu/\lambda$, $c_1=c/\lambda$.

From (2.2) $R(x,y,y,x) = R(x,Jy,Jy,x)$ and consequently
$$
	R\left( \frac{x+y}{\sqrt 2},\frac{x-y}{\sqrt 2},\frac{x-y}{\sqrt 2},\frac{x+y}{\sqrt 2} \right)=
	R\left( \frac{x+y}{\sqrt 2},\frac{Jx-Jy}{\sqrt 2},\frac{Jx-Jy}{\sqrt 2},\frac{x+y}{\sqrt 2} \right)
$$
which gives
$$
	H(x)+H(y)=4R(x,y,y,x)-2R(x,Jy,Jy,x)+2R(x,Jx,Jy,y)+2R(x,Jy,Jx,y) \ .
$$

Hence it is easy to find
$$
	(n+2)H(x) + \sum_{i=1}^{n} H(e_i) = S(x,x) +3S'(x,x)  \leqno(2.4)
$$
and
$$
	\sum_{i=1}^{n} H(e_i) =\frac{\tau(R)(p) +3\tau'(R)(p)}{4(n+1)}  \, . \leqno (2.5)
$$

From (2.4) and (2.5) we obtain
$$
	H(x) -\frac{1}{n+2}(S(x,x)+3S'(x,x))=-\frac{\tau(R)(p) +3\tau'(R)(p)}{4(n+1)(n+2)} \, .\leqno (2.6)
$$

Using (2.3) and (2.6) we get
$$
	\begin{array}{c}
		\ds \left( 2(n-2)\mu_1 - \frac{n+1}{n+2} \right) S(x,x)
		+\left( 2(n-1)\nu_1-\frac{3}{n+2} \right)S'(x,x)   \\
		=2(n-1)c_1(p) - \mu_1\tau(R)(p) -\nu_1\tau'(R)(p) -
		\ds \frac{\tau(R)(p) + 3\tau'(R)(p)}{4(n+1)(n+2)} \ .
	\end{array}            \leqno (2.7)
$$

Hence by a simple calculation we obtain
$$
	c_1 = \left( \frac{\mu_1}{n} + \frac{2n+1}{8n(n^2-1)} \right) \tau(R) 
	     +\left( \frac{\nu_1}{n} - \frac{3}{8n(n^2-1)} \right) \tau'(R) \ . \leqno (2.8)
$$

The substitution of (2.8) in (2.7) gives
$$
	\begin{array}{c}
		\ds \left( \mu_1 + \frac{n+1}{2(n^2-4)} \right) S(x,x)
		+\left( \nu_1-\frac{3}{2(n^2-4)} \right)S'(x,x)   \\
		=	\ds \frac{1}{2n}\left\{\left( \mu_1 +\frac{n+1}{2(n^2-4)}\right) \tau(R)(p) + 
		\left( \nu_1 - \frac{3}{2(n^2-4)} \right) \tau'(R)(p)\right\}  \ . 
	\end{array}                        \leqno (2.9)
$$

From (2.2), (2.8) and (2.9) it follows
$$
	\begin{array}{c}
		\ds R(x,y,y,x) -\frac{n+1}{2(n^2-4)}(S(x,x)+S(y,y))
		    +\frac{3}{2(n^2-4)} (S'(x,x) + S'(y,y))   \\
		\ds =-\frac{2n^2+3n+4}{8(n^2-1)(n^2-4)}\tau(R)(p) +
		    \frac{9n}{8(n^2-1)(n^2-4)} \tau'(R)(p) \ .  
	\end{array}                  \leqno (2.10)    
$$

According to the lemma, (2.6) and (2.10) imply that the Bochner curvature
tensor $B$ for $M$ vanishes. The rest of the theorem follows from (2.2) and (2.10).

\vspace {0.5cm}
\noindent
{\it Centre for Mathematics and Mechanics \ \ \ \ \ \ \ \ \ \ \ \ \ \ \ \ \ \ \ \ \ \ \ \ \ \ \ \ \ \ \ \ \ \
Received 26.5.1981

\noindent
Sofia 1090   \ \ \ \ \ \ \ \ \ \ \ \ \ \ \ \ \  P. O. Box 373}


\vspace{0.6in}

\begin{thebibliography}{9}

\bibitem{G}
G. G a n \v{c} e v. Almost Hermitian manifolds similar to the complex space forms. 
{\it C. R. Acad. bulg. Sci.}, {\bf 32}, 1979, 1179-1182.

\bibitem{S} 
J. S c h o u t e n, D. S t r u i k. On some properties of general manifolds relating to
Einstein's theory of gravitation. {\it Amer. J. Math.}, {\bf 43}, 1921, 213-216.

\bibitem{T} 
S. T a c h i b a n a. On the Bochner curvature tensor. {\it Nat. Sc. Rep. Ochanomizu
Univ.}, {\bf 18}, 1967, 15-19. 

\bibitem{TV}
F. T r i c e r r i, \, L. V a n h e c k e. Curvature tensors on almost Hermitian manifolds.
\ {\it Trans. Amer. Math. Soc.}, {\bf 267}, 1981, 365-398.

\bibitem{V}
L. V a n h e c k e. Almost Hermitian manifolds with $J$-invariant Riemann curvature tensor.
{\it Rend. Sem. Matem. Torino}, {\bf 34}, 1975-76, 487-498.

\end{thebibliography}
\end{document}